\documentclass[11pt]{article}
\usepackage{amsmath,amssymb,amsthm,graphicx,mathtools,url,verbatim}
\usepackage[margin=1in]{geometry}

\newtheorem*{zl*}{Zolotarev's Lemma}
\newtheorem*{qr*}{Law of Quadratic Reciprocity}

\newcommand{\ZZ}{\mathbf Z}
\newcommand{\sgn}{\operatorname{sign}}

\theoremstyle{definition}

\newtheorem*{exercise*}{Exercise}

\newcommand{\LS}[2]{{\genfrac{(}{)}{}{}{#1}{#2}}}
\newcommand{\ZS}[2]{{\genfrac{[}{]}{}{}{#1}{#2}}}

\title{Zolotarev's Magical Proof of Quadratic Reciprocity\thanks{The author was partially supported by NSF grant DMS-2154224.}}

\author{Matthew Baker}
\date{\today}

\begin{document}

\maketitle

\begin{quote}
Every great magic trick consists of three parts or acts. The first part is called ``The Pledge''. The magician shows you something ordinary\ldots Perhaps he asks you to inspect it to see if it is indeed real, unaltered, normal. But of course... it probably isn't. The second act is called ``The Turn''. The magician takes the ordinary something and makes it do something extraordinary\ldots But you wouldn't clap yet. Because making something disappear isn't enough; you have to bring it back. That's why every magic trick has a third act, the hardest part, the part we call ``The Prestige''.

--- Cutter (Michael Caine) in Christopher Nolan's ``The Prestige''

\end{quote}

\section{The Pledge}

Let $m$ and $n$ be positive integers.
Suppose you have a deck consisting of $mn$ different cards numbered $0,1,\ldots,{mn-1}$ from top to bottom.
Consider the following two ways of dealing the cards into an $m \times n$ rectangular grid:

\medskip

\medskip
\noindent
\textbf{(R)}: \emph{Row deal.} 
Deal the cards by rows, dealing the $m$ rows from top to bottom and the $n$ cards within each row from left to right.

\begin{figure}[htbp]
    \centering
    \includegraphics[scale=0.075]{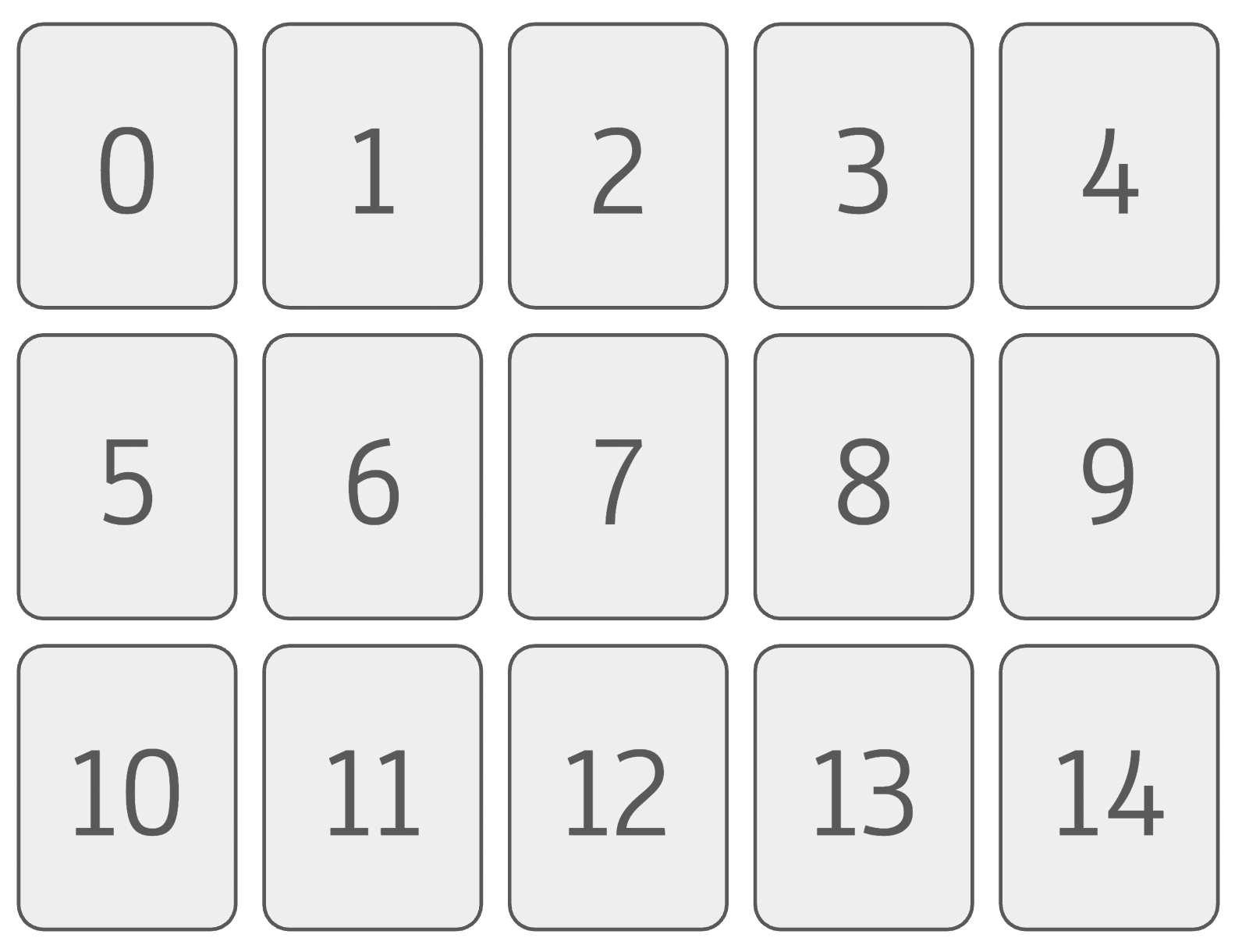}
    \caption{Example of a row deal with $m=3$ and $n=5$}
    \label{fig:Rho}
\end{figure}

\medskip
\noindent
\textbf{(C)}: \emph{Column deal.} 
Deal the cards by columns, dealing the $n$ columns from left to right and the $m$ cards within each column from top to bottom.

\begin{figure}[htbp]
    \centering
    \includegraphics[scale=0.1]{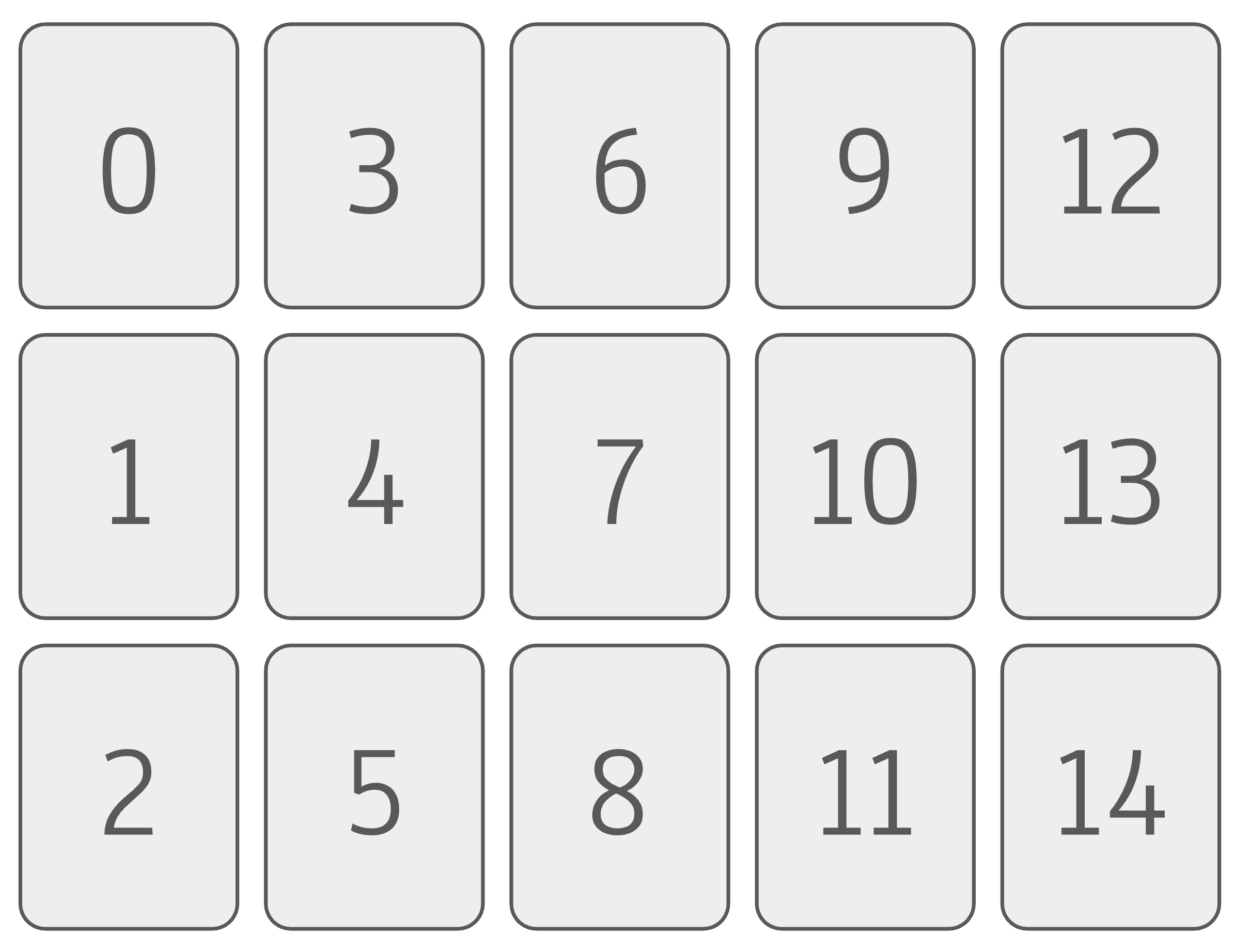}
    \caption{Example of a column deal with $m=3$ and $n=5$}
    \label{fig:Kappa}
\end{figure}

We label the rows of the grid by $0,1\ldots,m-1$ and the columns by $0,1,\ldots,n-1$.
The map which takes (R) to (C) yields a permutation $\gamma$ of the set $S = \{ 0,1,\ldots,mn-1 \}$ of cards by replacing the card in position $(i,j)$ of (R) with the card in position $(i,j)$ of (C).
(We can think of $\gamma$ as first picking up the cards by rows and then dealing them out by columns.)


Now for a combinatorial puzzle: what is the {\em sign} of the permutation $\gamma$?
(If you'd like to think about it before seeing the answer, stop reading now!)

\medskip

The answer is that 
\begin{equation} \label{eq:generalgammaformula}
\sgn(\gamma) = (-1)^{\binom{m}{2} \cdot \binom{n}{2}}.
\end{equation}

To prove \eqref{eq:generalgammaformula}, recall that the sign of a permutation $\sigma$ of a totally ordered finite set $X$
is equal to $(-1)^{\iota({\sigma})}$, where $\iota(\sigma)$ is the number of {\em inversions} of $\sigma$.
(An inversion is an ordered pair $(x,y) \in X^2$ with $x<y$ and $\sigma(x) > \sigma(y)$.)

\begin{exercise*}
Show that a pair $(x,y)$ of cards in (R) form an inversion of $\gamma$ if and only if $y$ lies below and to the left of $x$.
\end{exercise*}

The number of inversions of $\gamma$ is therefore $\binom{m}{2} \cdot \binom{n}{2}$,
since each pair consisting of a $2$-element subset $\{ i,i' \}$ of $\{0,1,\ldots,m-1\}$ and a 2-element subset $\{ j,j' \}$ of $\{0,1,\ldots,n-1\}$ gives rise
to a unique inversion $(x,x')$, where $x=(i,j)$ and $x'=(i',j')$, by ordering the elements so that $i < i'$ and $j > j'$.  

\medskip

If $m$ and $n$ are both odd, then \eqref{eq:generalgammaformula} simplifies: in this case, we obtain
\begin{equation} \label{eq:gammaformula}
\sgn(\gamma) = (-1)^{\frac{(m-1)(n-1)}{4}}.
\end{equation}

In other words, $\sgn(\gamma) = 1$ if either $m \equiv 1 \pmod{4}$ or $n \equiv 1 \pmod{4}$,
and $\sgn(\gamma) = -1$ if $m \equiv n \equiv 3 \pmod{4}$.

\medskip

Formula~\eqref{eq:gammaformula} may bring to mind Gauss's Law of Quadratic Reciprocity.  
Is this just a coincidence?  Continue on, dear reader\ldots

\section{The Turn}

We now assume that $m$ and $n$ are {\em relatively prime}, in addition to being odd and positive.
This gives us a third way of dealing the cards into an $m \times n$ rectangular grid:

\medskip

\medskip
\noindent
\textbf{(D)}: \emph{Diagonal deal.} 
Starting in the top left corner, deal the cards diagonally down and to the right, wrapping around from bottom to top or right to left whenever necessary.\footnote{For the topologists in the audience, we are actually dealing the cards onto a torus.}

\begin{figure}[htbp]
    \centering
    \includegraphics[scale=0.1]{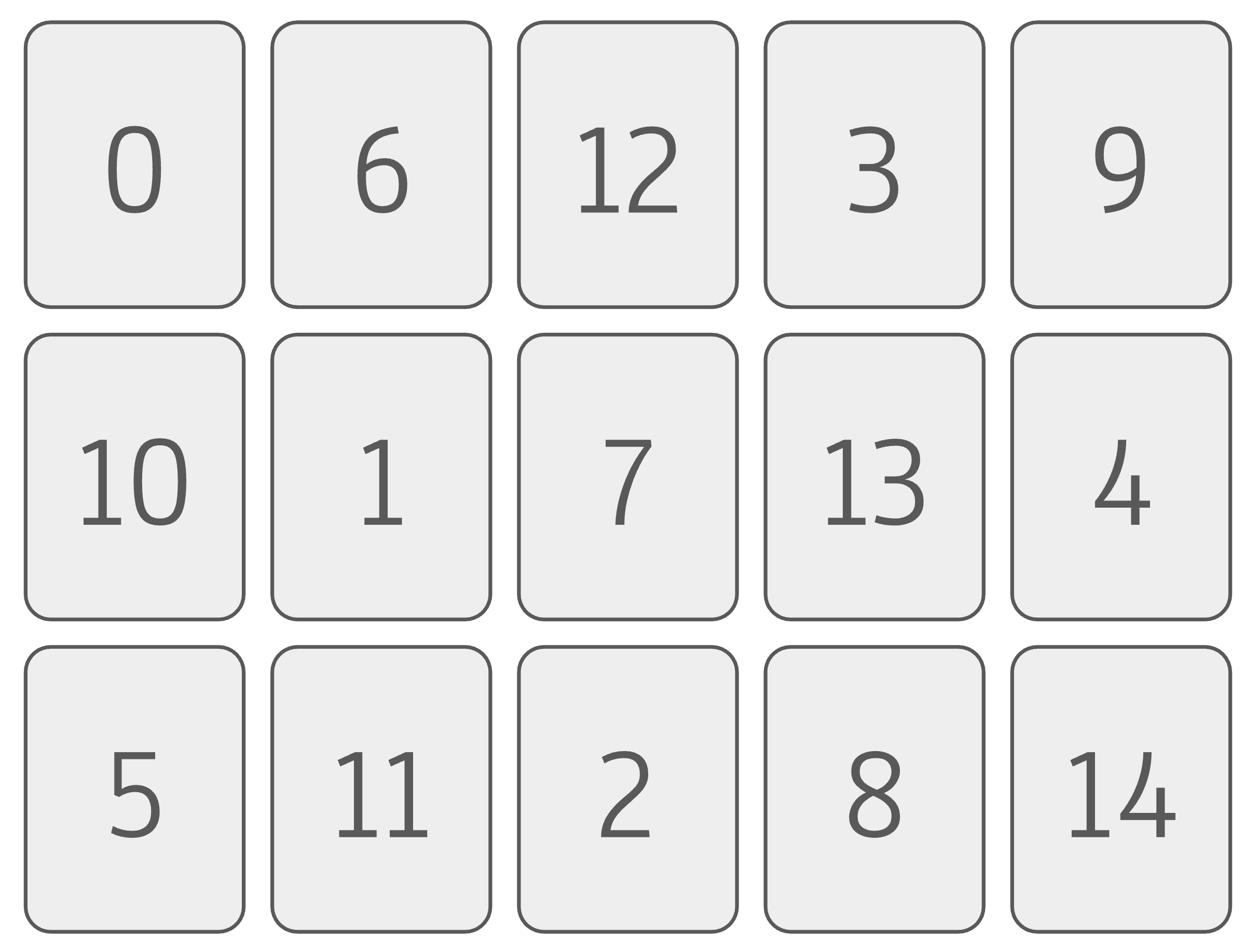}
    \caption{Example of a diagonal deal with $m=3$ and $n=5$}
    \label{fig:Delta}
\end{figure}

The map which takes (R) to (D) yields a permutation $\alpha$ of $S$.
(We can think of $\alpha$ as first picking up the cards by rows and then dealing them out diagonally.)

\medskip

Now for another combinatorial puzzle: what is the sign of $\alpha$?
(Again, if you'd like to think about it before seeing the answer, stop reading now!)

\medskip

We claim that $\sgn(\alpha)$ is equal to the sign of the permutation of $\ZZ / m \ZZ$ induced by multiplication by $n$.
In other words, writing $\ZS{n}{m}$
for the sign of this permutation, we are claiming that
\begin{equation} \label{eq:alphaformula}
\sgn(\alpha) = \ZS{n}{m}.
\end{equation}

To prove \eqref{eq:alphaformula}, we begin with the following observations:

\begin{exercise*}
\begin{enumerate}
\item[]
\item Show that the card in position $(i,j)$ of (D) is the unique card congruent to $i$ modulo $m$ and to $j$ modulo $n$.
\item Show that the card in position $(i,j)$ of (R) is $ni + j$.
\item Show that $\alpha$ permutes each individual column of (R), and $\alpha^{-1}$ restricted to column $j$ sends the card in row $i$ to the unique card congruent to $ni + j$ modulo $m$.
\end{enumerate}
\end{exercise*}

It follows that $\alpha^{-1}$ is a product of $n$ permutations, each of which can be expressed as a composition of multiplication by $n$ modulo $m$ and translation by $j$ modulo $m$.
Since $m$ is odd and translation by $j$ modulo $m$ is a product of $k$-cycles with $k$ dividing $m$, it follows that translation by $j$ modulo $m$ is an even permutation for every $j$.
Therefore, since the sign of a permutation is a homomorphism, we have $\sgn(\alpha) = \sgn(\alpha^{-1})=  \ZS{n}{m}^n$. Since $n$ is also odd, this establishes \eqref{eq:alphaformula}.

\medskip

We now swap the roles of rows and columns in the previous argument.
The map which takes (D) to (C) yields a permutation $\beta$ of $S$.
Since $\sgn(\beta) = \sgn(\beta^{-1})$ and $\beta^{-1}$ can be obtained from $\alpha$ by swapping the roles of $m$ and $n$, it follows by symmetry that
$\sgn(\beta) = \ZS{m}{n}$.
And since $\beta \circ \alpha = \gamma$, we have
\begin{equation*} \label{eq:signcomposition}
\sgn(\beta) \cdot \sgn(\alpha) = \sgn(\gamma),
\end{equation*}
i.e.,
\begin{equation} \label{eq:signcomposition2}
\ZS{m}{n} \cdot \ZS{n}{m} = (-1)^{\frac{(m-1)(n-1)}{4}}.
\end{equation}

\medskip

Formula~\eqref{eq:signcomposition2} is even \emph{more} reminiscent of the Law of Quadratic Reciprocity!
But wait\ldots  what is the connection to perfect squares modulo primes?

\medskip

Now that you are hooked, dear reader, you have no choice but to continue\ldots

\section{The Prestige}

Let $p$ be an odd prime, and let $a$ be a positive integer not divisible by $p$. Recall that the \emph{Legendre symbol} $\LS{a}{p}$
is defined by 
\[
\LS{a}{p} = 
\begin{cases}
1 & \text{if $a$ is a square modulo $p$} \\
-1 & \text{if $a$ is not a square modulo $p$}.\\
\end{cases}
\]

The connection between the ``Zolotarev symbol'' $\ZS{a}{p}$ and the Legendre symbol $\LS{a}{p}$ is given by:

\begin{zl*}
If $p$ is an odd prime and $a$ is a positive integer not divisible by $p$, then
\[
\ZS{a}{p} = \LS{a}{p}.
\]
\end{zl*}

To prove Zolotarev's Lemma, we first recall Gauss's \emph{primitive root theorem}: if $p$ is an odd prime, the multiplicative group $(\ZZ / p\ZZ)^\times$ is \emph{cyclic}.
It follows that $\ZS{\cdot}{p}$ is a surjective homomorphism from the multiplicative group $(\ZZ / p\ZZ)^\times$ to $\{ \pm 1 \}$;
surjectivity follows from the fact that if $g$ is a primitive root mod $p$ (i.e., a cyclic generator of $(\ZZ / p\ZZ)^\times$) then 
$\ZS{g}{p}$ is a $(p-1)$-cycle and thus has signature $-1$. 
The kernel of $\ZS{\cdot}{p}$ is therefore a subgroup of  $(\ZZ / p\ZZ)^\times$ of index $2$,
but the only such subgroup is the group of quadratic residues.  Thus
$\ZS{\cdot}{p}$ coincides with the Legendre symbol $\LS{\cdot}{p}$.

\medskip

Combining Zolotarev's Lemma with \eqref{eq:signcomposition2} yields:

\begin{qr*}
\label{cor:quadrecip}
If $p$ and $q$ are distinct odd primes, then
\[
\boxed{\displaystyle
\LS{p}{q} \cdot \LS{q}{p} = (-1)^{\frac{(p-1)(q-1)}{4}}.
}
\]
\end{qr*}

\section*{Closing Credits}

Gauss found eight different proofs of the Law of Quadratic Reciprocity over the course of his life. 
There are (at the present moment) 345 different proofs listed in Franz Lemmermeyer's database \cite{lemmermeyerQRG}.
The proof I've given here is essentially a creative reimagining and dealgebraization of a beautiful 1872 proof due to Egor Ivanovich Zolotarev (1847--1878) \cite{zolotarev1872}, which is not as well-known as it ought to be. My ``card dealing'' 
interpretation of Zolotarev's proof first appeared on my blog in 2013 \cite{baker2013zolotarev}.
Thanks to Timothy Chow, Jerry Shurman, and Joe Silverman for their feedback on the exposition.

\section*{Director's Commentary}
\begin{enumerate}
\item Zolotarev defended his doctoral thesis in 1874 and was appointed as an ``extraordinary professor'' by the Russian Academy of Sciences in 1876, but his meteoric rise in the math world ended abruptly when he was run over by a train at the age of 31. 
 \item If we identify the set $\{0,1,2,\ldots,mn-1\}$ with the set $\ZZ/mn\ZZ$ of integers modulo $mn$, the map $\gamma$ is just the canonical ring isomorphism from $\ZZ/mn\ZZ$ to $\ZZ/m\ZZ\times \ZZ/n\ZZ$ afforded by the Chinese Remainder Theorem.
 \item Zolotarev's Lemma generalizes to the statement that if $m,n$ are relatively prime odd positive integers then $\ZS{m}{n}$ is equal to the Jacobi symbol $\LS{m}{n}$, and the proof above then gives quadratic reciprocity for the Jacobi symbol (which is used for rapid computation of Legendre symbols). 
\end{enumerate}

\section*{Director's Cut Bonus Footage}

We show how to use similar arguments to obtain the two ``supplements'' to the Law of Quadratic Reciprocity, which give formulas for $\ZS{2}{n}$ and $\ZS{-1}{n}$, respectively.

Let $n \ge 3$ be an odd integer and put $m = \frac{n-1}{2}$.  We now work with a deck of
$n-1$ cards labeled $1,2,\dots,n-1$, which we deal into a $2\times m$ rectangular array in three different ways as follows:

\medskip
\noindent
\textbf{(R)}: \emph{Row deal.} Deal the cards by rows, as before.

\smallskip
\noindent
\textbf{(C)}: \emph{Column deal.} Deal the cards by columns, as before.

\smallskip
\noindent
\textbf{(Z)}: \emph{Zigzag deal.} Deal the first card into the lower-left corner, then alternate between moving up and moving down and to the right.
Observe that (Z) can be obtained from (C) by simply swapping the two entries in each column.

\medskip

Let $\gamma : (R) \to (C)$ as before, and let $\alpha : (R) \to (Z)$, $\beta : (Z) \to (C)$.
Then $\gamma = \beta\circ\alpha$, hence
\[
\operatorname{sign}(\gamma)=\operatorname{sign}(\beta)\cdot\operatorname{sign}(\alpha).
\]

By our previous calculation \eqref{eq:generalgammaformula}, we have
\[
\operatorname{sign}(\gamma)=(-1)^{\binom{2}{2} \binom{m}{2}} = (-1)^{\binom{m}{2}}.
\]

Since $(Z)$ is obtained from $(C)$ by swapping the entries in each of the $m$ columns,
$\beta$ is a product of $m$ transpositions and thus
\[
\operatorname{sign}(\beta)=(-1)^{m}.
\]

By direct inspection, we see that $\alpha$ corresponds
to multiplication by $2$ modulo $n$ on the set $\{1,2,\dots,n-1\}$.
Thus $\operatorname{sign}(\alpha)=\ZS{2}{n}$.

Combining the above identities gives
\[
\ZS{2}{n}
=\operatorname{sign}(\alpha)
=\operatorname{sign}(\beta^{-1})\operatorname{sign}(\gamma)
=\operatorname{sign}(\beta)\operatorname{sign}(\gamma)
=(-1)^{m+\binom{m}{2}}
=(-1)^{\frac{m(m+1)}{2}}
=(-1)^{\frac{n^{2}-1}{8}}.
\]

By Zolotarev’s lemma, if $n=p$ is an odd prime then $\ZS{2}{p}=\LS{2}{p}$,
so we recover the classical supplemental law:
\[
\boxed{\displaystyle
\LS{2}{p} = (-1)^{\frac{p^2-1}{8}}.}
\]

\medskip

\begin{exercise*}
With $n$ and $m$ as above, we introduce one more method of dealing the cards into a $2\times m$ rectangular array:

\smallskip
\noindent
\textbf{(M)}: \emph{Modified zigzag deal.} Deal the first card into the lower-right corner, then alternate between moving up and moving down and to the left.
\begin{enumerate}
\item Show that the permutation $\delta : (R) \to (M)$ corresponds 
to multiplication by $-1$ modulo $n$ on the set $\{1,2,\dots,n-1\}$, hence $\operatorname{sign}(\delta)=\ZS{-1}{n}$.
\item Show that (M) can be obtained from (R) by first reversing the order of the columns, then swapping the two entries in each column. Conclude that $\operatorname{sign}(\delta)=(-1)^m$.
\item Using Zolotarev's lemma, conclude that if $n=p$ is an odd prime then $\LS{-1}{p} =(-1)^{\frac{p-1}{2}}$.
\end{enumerate}
\end{exercise*}

\end{document}